\documentclass[12pt,a4paper]{article}

\usepackage[dvips]{graphicx}
\usepackage{floatflt}
\usepackage{wrapfig}
\usepackage{longtable}

\def\bb{\begin{equation}}
\def\ee{\end{equation}}

\def\const{\hbox{const}}

\def\ve{\varepsilon}

\newtheorem{theorem}{Theorem}
\newtheorem{lemma}{Lemma}
\newtheorem{corollary}{Corollary}
\title{Oscillations near separatrix for perturbed Duffing equation}

\author{O.M. Kiselev,
\footnote{Institute of Mathematics 
USC of RAS, ok@ufanet.ru}}
\date{\null}

\begin{document}
\maketitle
\begin{abstract}
A periodic perturbation generates a
complicated dynamics close to separatrices and saddle points.  We construct 
an asymptotic solution which is close to the separatrix  for
the unperturbed Duffing's oscillator over a  long time. This solution is defined by 
a separatrix map. This map is obtained for any order of the perturbation parameter.
Properties of this map show an instability of a motion for the perturbed system.
\end{abstract}
\par
{\Large \bf Introduction}
\par
We consider an equation for the perturbed Duffing's oscillator:
\begin{eqnarray}
u''+2u-2u^3=\ve\cos(\omega t+\Phi_0).
\label{pertubedDuffingsOscillator}
\end{eqnarray}
Here $\ve$ is a small positive parameter, $\omega$ and $\Phi_0$ are constants. 
\par
The goal of the paper is to construct an asymptotic solution for
(\ref{pertubedDuffingsOscillator}) which is close to the separatrix  for
the unperturbed Duffing's oscillator over a  long time. Formally the words 
{\it asymptotic solution} mean that there exists $\ve_0=\const>0$, such that  
a following asymptotic series 
\begin{equation}
U(t,\ve)=\sum_{n=0}^{\infty}\ve^n 
U_n(t),
\label{separatrixAsymptoticExpansion}
\end{equation}
gives a residual as $o(\ve^n), \quad \forall n\in{\mathbf N},$
as $\ve\in (0,\ve_0)$ when $U(t,\varepsilon)$ is being substituted into
(\ref{pertubedDuffingsOscillator}). 
\par
Our plane is to concentrate on asymptotic solution with a separatrix as
leading term:
\begin{equation}
U_0(t)=\tanh(t+\tilde t).
\label{upperSeparatrix}
\end{equation}
\par
A similar problem was studied by A. Poincare
\cite{Poincare}. He considered a special case of  the three body
problem. Later V.K. Mel'nikov established that the separatrices split under a 
perturbation and he calculated a gap between the separatrices of
perturbed equation \cite{Melnikov}. N.N.Filonenko, R.Z.Sagdeev and G.M. Zaslavskii
\cite{FilonenkoSagdeevZaslavskii} obtained a separatrix map for
canonical variables in Hamiltonian systems. 
\par
The problem for the perturbation of a separatrix is very important for a 
capture into a resonance. Therefore the problem of dynamics near the
separatrix was studied in the point of view of the
capture. A.I. Neishtadt \cite{Neishtadt} calculated a measure of
trajectories which cross the separatrix and are captured into a 
resonance.  A.V.Timofeev studied the dynamics near the separatrix
close to the saddle  \cite{Timofeev}. 
\par
A change of the angle variable when the solution crosses the
separatrix was obtained by J.R.Cary and R.T.Scodjie
\cite{CaryScodjie}. D.C.Diminie and R.Haberman studied a separatrix
crossing near a saddle-center bifurcation and a pitchfork bifurcation
\cite{DiminieHaberman1}. Full asymptotic expansions
for the problem of  the separatrix crossing near the saddle-center bifurcation
was obtained in \cite{Kiselev,KiselevGlebov}.
\par
In this work  following new results are obtained. It is constructed
the map for the dynamics near the separatrix for any power of the perturbation
parameter.   It is shown that the phase shift over a circle near
the separatrix circle is defined by a term of the order  $\ve^2$ on the previous circle. The
constructed map shows that  there exists a manifold of
$\ve\in(0,\ve_0)$ for $\forall \ve_0>0$, such that the asymptotic solution has more than
$N$ circles near the separatrix for $\forall N\in{\mathbf N}$. A Cantor set gives
an example of such manifolds. 
\par
This work has following structure. In section
\ref{mainProblemAndResults} the main problem and results are presented
in formal form. Section \ref{separatrixDynamics} contains the derivation of the map for  the dynamics near
the separatrix. In section \ref{sectionDescreetDynamicalSystem}
the consequences  of the separatrix 
dynamics are presented. 

\section{Main problem and results}
\label{mainProblemAndResults}
\subsection{A trial solution}
\par
Let us consider an asymptotic solution in the form
(\ref{separatrixAsymptoticExpansion}). Main problem is to construct a
bounded asymptotic solution of such type as $t\in(0,-N\ln(\ve))$ for
any $N\in{\mathbf N}$, and $\varepsilon\in(0,\varepsilon_0)$. 
\par
A simplest result is following.
\begin{theorem}
\label{theoremAboutSeparatrixPerturbation}
There exists two parametric asymptotic  solution of
(\ref{pertubedDuffingsOscillator}) in form 
(\ref{separatrixAsymptoticExpansion}) and (\ref{upperSeparatrix}) when 
$$
{1\over2}\ln(\ve)\ll t\ll -{1\over2}\ln(\ve).
$$
Higher-order terms of (\ref{separatrixAsymptoticExpansion}) are 
$$
U_n(t)=A_n^- e^{2t}+B_n^- e^{-2t}+W^-_n(t),
$$
as $t\to-\infty$, where
$$
W_n^-(t)=O(e^{-2nt}),\quad t\to-\infty,
$$
and $W_n^-(t)$ has an asymptotic expansion into series of powers  of
$e^t$ as $t\to\infty$, such that it does not contains the terms  $C_1e^{2t}$
and $C_2e^{-2t}$  where $C_1$ and $C_2$ are some constants. This means
that all parameters for $U_n(t)$ are $A_n^-$ and $B_n^-$. There fore
the parameters of the asymptotic solution are:
\begin{eqnarray*}
A^-=\sum_{n=1}^\infty \varepsilon^n A_n^-,\\
B^-=\sum_{n=1}^\infty \varepsilon^n B_n^-.
\end{eqnarray*}
\end{theorem}
\par
Parameter $A^-$ can be excluded from the solution by time shift in
main term of the asymptotic solution (\ref{upperSeparatrix}) and
$\Phi_0$ in the  perturbation term of
(\ref{pertubedDuffingsOscillator}). Therefore we take a case
$A_n^-=0$. Parameter $B^-$ is a distance between the asymptotic
solution and the separatrix of unperturbed equation (\ref{pertubedDuffingsOscillator}). 

\subsection{Numeric simulations and challenge for the analytic studies}
\par
\begin{wrapfigure}{r}{0.5\textwidth}
\vspace{-30pt}
\includegraphics[width=8cm,height=6cm]{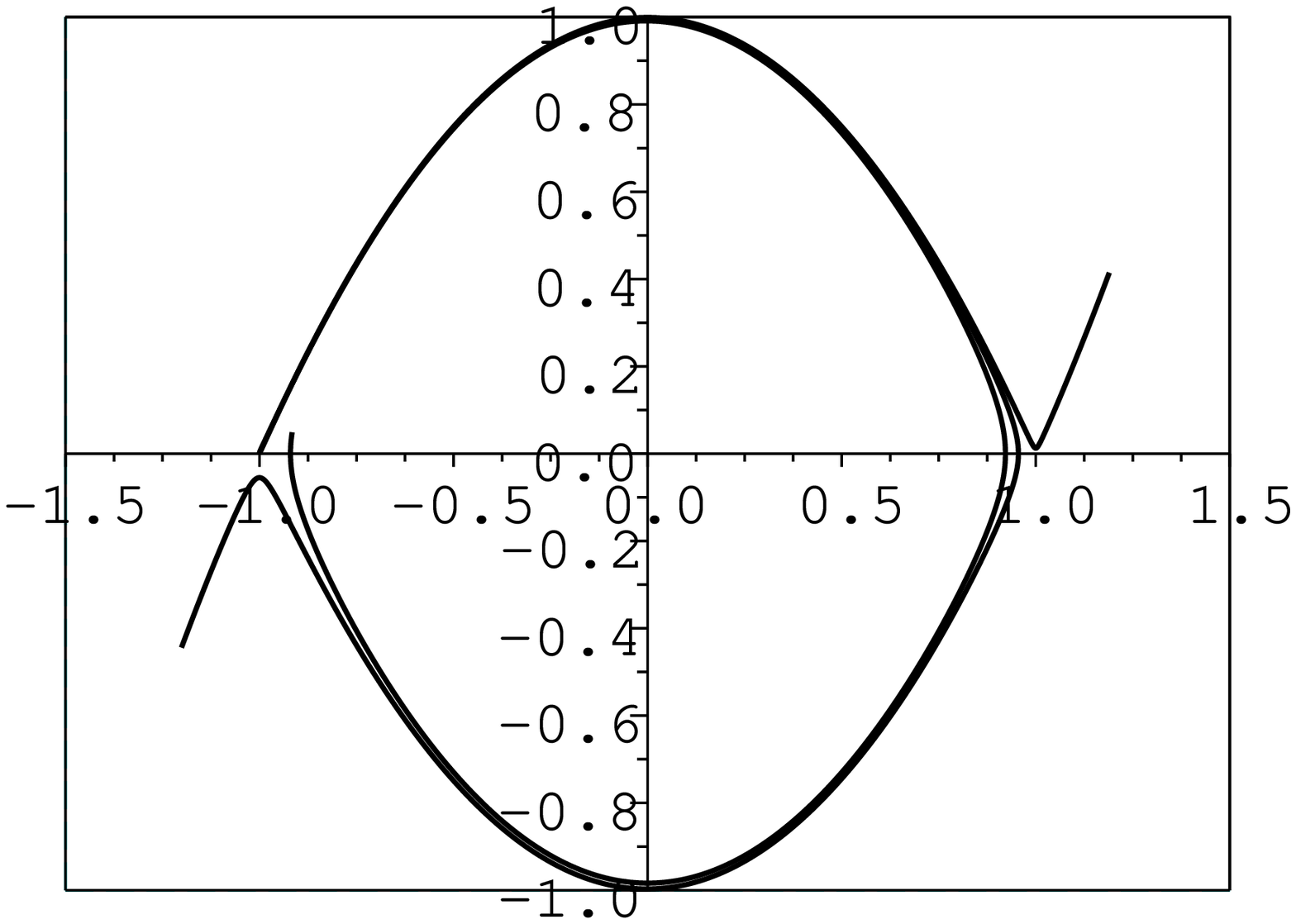}
\vspace{-30pt}
\caption{\small Trajectories near separatrix on the plane $(u,u')$.}\label{neighbourOfRightSaddle}
\vspace{-20pt}
\end{wrapfigure}
There are two scenarios for a prolongation of the trial solution. When
$t=O(-\ln(\varepsilon)/2)$ the trajectory of the trial solution closes to
the saddle point $(1,0)$. This trajectory is able to turn to the lower
separatrix, which goes to another saddle point
$(-1,0)$. The different way for the trajectory lies near the separatrix which
goes from $(1,0)$ to $(+\infty,+\infty)$.  The similar changes are possible near 
the saddle $(-1,0)$. These cases are shown on the figure \ref{neighbourOfRightSaddle}.
\par
\begin{figure}
\hspace{-1cm}
\includegraphics[width=0.55\textwidth]{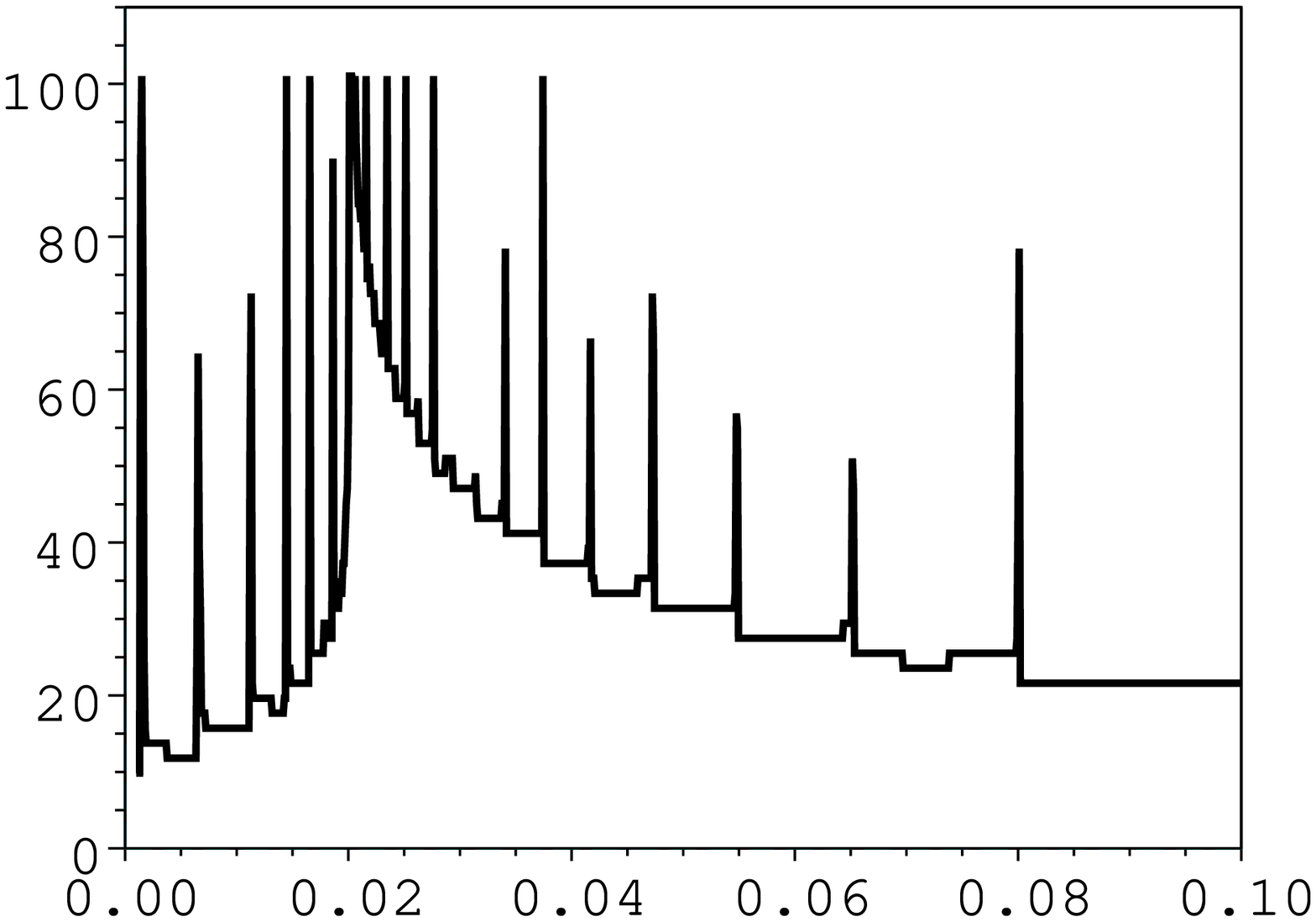}
\includegraphics[width=0.55\textwidth]{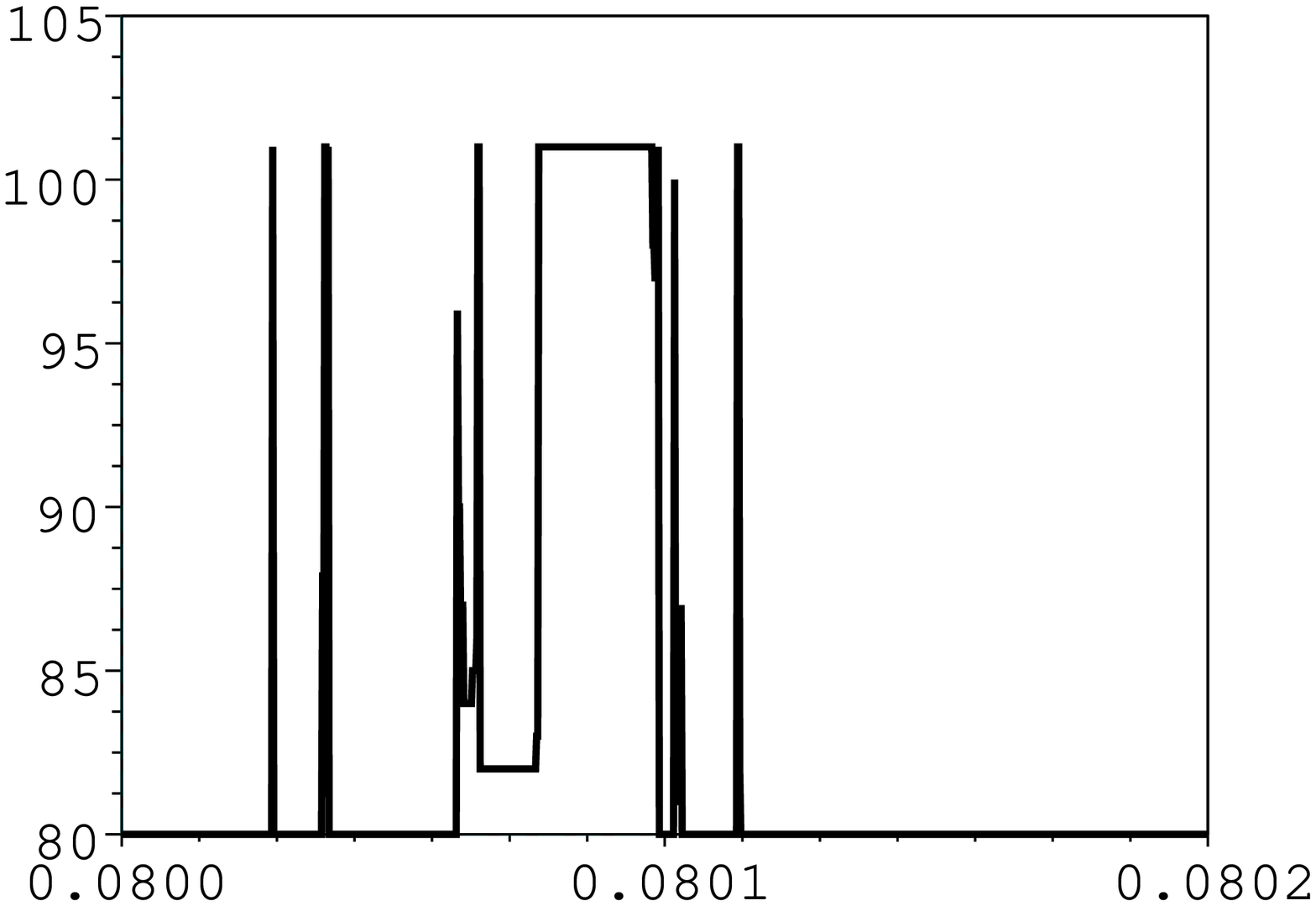}
\caption{\small The left figure shows the dependence of the oscillating time on $\varepsilon$  for the solution of the Cauchy problem $u(0)=0,\, u'(0)=1$ for (\ref{pertubedDuffingsOscillator}) where $\omega=1,\,\Phi_0=\pi$. The solutions was studied on the interval $t\in(0,100)$ for the values of the perturbation parameter $\varepsilon\in(0.001,0.1)$. The right figure shows the scaled structure of the peak near $\varepsilon=0.08$. \label{ocsillationsTime1}}
\end{figure}

Let us concentrate to the trajectories which oscillates between the
saddles $(-1,0)$ and $(1,0)$. For this we should to find manifolds
for the parameters of a solution and the perturbation $\varepsilon$. These
manifolds have a complicated structure. For example one can see the
dependency of the life-time for oscillating asymptotic solution on the
perturbation parameter $\varepsilon$. This shown on the right figure
\ref{ocsillationsTime1}. On left picture
\ref{ocsillationsTime1} one can see the thin structure of the peak
near $\varepsilon=0.08$.
\par
This numeric simulations take a challenge for analytic studies. Our
goals are to calculate the asymptotic solutions, to find a dependency
for the trajectory on the parameters and to give the formulas for  
the manifolds of the parameters for the solutions with an oscillating behaviour.

\subsection{Results}
\par
The asymptotic solution which is defined in Theorem
\ref{theoremAboutSeparatrixPerturbation} will be used as a trial
solution and we  will prolong it to the large time. The asymptotic
behaviour of this solution has the same form:
$$
U_n(t)=A_n^+ e^{-2t}+B_n^+ e^{2t}+W^+_n(t),
$$
as $t\to\infty$, where
$$
W_n^+(t)=O(e^{-2nt}),\quad t\to+\infty,
$$
We will show if $B^+_1<0,$ then the solution (\ref{separatrixAsymptoticExpansion}),
(\ref{upperSeparatrix}) is bounded when 
$$
{1\over2}\ln(\varepsilon)\ll t\ll -\ln(\ve)
$$
and when $1\ll t\ll
-\ln(\varepsilon)$ the solution has a following asymptotic expansion:
\begin{eqnarray*}
u(t,\ve)=-\tanh(\theta)+
\sum_{n=1}^\infty \ve^n u_n(\theta),\\
\theta=t+{1\over2}\ln(\ve)+
{1\over2}\ln\big(-{1\over16} B_1^+\big),
\end{eqnarray*}
as $\theta\to\mp\infty$, where 
$$
u_n(\theta)=a_n^\mp e^{\pm2\theta}+
b_n^\mp e^{\mp2\theta}+w^\mp_n(\theta),
$$
where $w_n^\mp(\theta)$ has the same properties as $W_n^\pm(t)$ in Theorem \ref{theoremAboutSeparatrixPerturbation} 
\par
The following theorem gives conditions for the prolongation 
of previous theorem on $N$-circles near the separatrix 
of the unperturbed equation. 
\begin{theorem}
\label{sequencesOfSeparatrixCircle}
The asymptotic solution of (\ref{pertubedDuffingsOscillator}) from Theorem
\ref{theoremAboutSeparatrixPerturbation} may be extended on the interval $-k\ln(\ve)\ll t\ll 
-(k+{1\over2})\ln(\ve)$, $k=0,1,\dots,\big[{N\over2}\big]$ and has a form
$$
u(t,\ve)=(-1)^{k-1}\tanh(t_k)+
\sum_{n=0}^\infty \ve^n u_n^k(t_k),
$$
where $u_n^k(t_k)$ has the following asymptotic behaviour 
$$
u_n^k(t_k)\sim \sum_{\kappa=-\infty}^{n} e^{2\kappa t_k}\bigg(\sum_{l=0}^{2n}
  t_k^l\bigg[\sum_{m=0}^{n+1}\bigg( (\tilde u_n^k)^{\pm}_{\kappa,l,m}\cos(m\omega
  t_k)+ (\tilde v_n^k)^\pm_{\kappa,l,m}\sin(m\omega
  t_k)\bigg)\bigg]\bigg),
$$
as $t_k\to\pm\infty$. The condition of extendability is $(-1)^m\big[\sigma_1(m)+\Delta\sigma_1(m)\big]<0$ for $\forall
m=1,\dots,N,\,\, N\in{\mathbf N}$. Here the parameter
$\Delta\sigma_1(k)$ is following
\begin{eqnarray*}
\Delta\sigma_1(k)={\pi\over16\cosh(\pi \omega/2)}
\cos(\psi_k),
\end{eqnarray*}
and $t_k$, $\psi(k)$ and $\sigma_1(k)$ are defined by the recurrent sequence:
\begin{eqnarray*}
t_1=t,\quad t_{k+1}=t_k+{\omega\over2}\ln(\varepsilon)+\psi(k),\\
\sigma_k(1)=B^-_k,\quad \psi(1)=\phi,\quad \chi_k(1)=A^-_k,\\
\sigma_k(n+1)=-32{\sigma_{k+1}(n)+\Delta\sigma_{k+1}(n)\over
\sigma_{1}(n)+\Delta\sigma_{1}(n)},\\
\psi(n+1)=-{\omega\over2}\bigg(\ln\big({1\over16}(\sigma_{1}(n)+\Delta\sigma_{1}(n))\big)-\ln(2)\bigg)+\psi(n),\\
\chi_1(n)=-2,\quad 
\chi_k(n+1)=-{1\over32}(\sigma_{1}(n)+\Delta\sigma_{1}(n))(\chi_{k-1}(n)+\Delta\chi_{k-1}(n)),\\
\end{eqnarray*}
Here 
\begin{eqnarray*}
\Delta\sigma_k(n)=(u_n^k)^+_{1,0,0}-(u_n^k)^-_{-1,0,0},\\
\Delta\chi_k(n)=(u_n^k)^+_{-1,0,0}-(u_n^k)^-_{1,0,0}.
\end{eqnarray*}
\end{theorem}
\par
\begin{corollary}
The main term of the asymptotic expansion for $u(t,\varepsilon)$ depends
on $B_j^-,\ k=1,\dots,j$ as $t_k=O(1)$. This means that the asymptotic
expansion is unstable with respect to a small correction of the parameters.
\end{corollary}
\begin{corollary}
There exist a set of the parameters $\{A^-_k\}_{k=1}^\infty$,
$\{B^-_k\}_{k=1}^\infty $ such that the solution has an oscillating behaviour near
 the separatrices as $\delta=\omega\ln(\varepsilon)$ belongs by a Cantor set.
\end{corollary}

\section{Separatrix dynamics}
\label{separatrixDynamics}
\par
In this section theorem \ref{theoremAboutSeparatrixPerturbation} is
proved in three steps. First step is a construction of the asymptotic
expansion which is valid near the saddle points $u=\pm1$. On the second
step we obtain an asymptotic expansion which is valid close to
the separatrices $u=\pm\tanh(t)$. On last step we match these asymptotic
expansions and construct the uniform asymptotic solution, which  is
valid over all domains close to the separatrices and the saddles.

\subsection{Separatrix asymptotic 
expansion}
\label{subsectionSeparatrixAsymptoticExpansion}
\par
Let us to construct an asymptotic solution in the form:
\begin{equation}
u(t,\ve)=\sum_{n=0}^\infty \ve^n 
U_n(t).
\label{separtrixExpansion}
\end{equation}
\par
The leading term of the expansion is a separatrix solution 
$$
U_0=\tanh(t)
$$
of the unperturbed equation:
$$
U_0''+2U_0-2U_0^3=0.
$$
We will study the asymptotic expansion where the main term is the separatrix solution of 
\begin{equation}
U_n''+2 U_n-6 U_0^2 U_n=f_n,
\label{linearizedEq}
\end{equation}
where
$$
f_1=\cos(\omega t+\phi_0),
$$
and  when $n>1$  $f_n$ is a polynomial of third order with respect to $u_j$, $u_k$ and $u_l$ for $j+k+l=n$. 
\par
There are two linear independent solutions of the linearized equation:
$$
v''+2v-6 U_0^2 v=0.
$$
There solutions are:
$$
v_1={1\over \cosh^2(t)},\quad 
v_2={\sinh(4t)\over 32\cosh^2(t)} + 
{\sinh(2t)\over 4\cosh^2(t)} +
{3 t\over 8\cosh^2(t)}.
$$
Wronskian of the functions $v_1$ and $v_2$ is equal to unit. 
\par
A general solution of (\ref{linearizedEq}) has a following form:
\begin{equation}
U_n=
v_1(t)\int_{t_0}^t d\tilde t f_n(t) v_2(\tilde t) - 
v_2(t)\int_{t_0}^t d\tilde t f_n(t) v_1(\tilde t) +A_1 v_1 +B_1 v_2.
\label{solutionN-th}
\end{equation}
Here $t_0$ is a constant and $A_1$, $B_1$ are parameters of the
solution. 
\begin{lemma}
Let 
$$
f_n\sim \sum_{k=-\infty}^{n-1} e^{2kt}\bigg(\sum_{l=0}^{2(n-1)}
  t^l\bigg[\sum_{m=0}^n\bigg(F^{\pm}_{k,l,m}\cos(m\omega
  t)+H^\pm_{k,l,m}\sin(m\omega t)\bigg)\bigg]\bigg),\quad
  t\to\pm\infty,
$$
then 
$$
U_n(t)\sim \sum_{k=-\infty}^{n} e^{2kt}\bigg(\sum_{l=0}^{2n}
  t^l\bigg[\sum_{m=0}^{n+1}\bigg(\tilde U^{\pm}_{k,l,m}\cos(m\omega
  t)+\tilde V^\pm_{k,l,m}\sin(m\omega t)\bigg)\bigg]\bigg),\quad t\to\pm\infty.
$$
\end{lemma}
To prove this lemma one should use asymptotic behaviours for
$v_{1,2}$, substitute the formula for the $f_n$
into (\ref{solutionN-th}) and integrate. 
\par
Denote $A_n^{\pm}=\tilde U_{\mp1,0,0}^{\pm}/4$ and $B_n^{\pm}=16\tilde
U_{\pm1,0,0}^{\pm}$. These parameters define the solution for a large
time. The changes of these parameters are:
$$
\Delta A_n=A_n^+ -A_n^-,\quad \Delta B_n=B_n^+ - B_n^-.
$$
The value of $\Delta B_1$ is defined by Melnikov's integral:
\begin{eqnarray*}
\Delta B_1 &=B_1^+ -B_1^- = \int_{-\infty}^{\infty}
{\cos(\omega t+\Phi_0)\over\cosh^2(t)}dt
\\
&=\int_{-\infty}^{\infty}
{dt\over \cosh^2(t)}
(\cos(\omega t)\cos(\Phi_0)-
\\
&-\sin(\omega t)\sin(\Phi_0))
=\cos(\Phi_0){\pi\over
\cosh(\pi\omega/2)}.
\end{eqnarray*}
The following formula gives $\Delta A_1$:
\begin{eqnarray*}
\Delta A_1=\lim_{s\to+\infty}\bigg[-\sin(\Phi_0)\bigg(\int_{-s}^{s}v_2(t)\sin(\omega t)dt-\\
{2\sin(\omega s)-\omega\cos(\omega s)\over 8(\omega^2+4)}e^{2s}
\bigg)\bigg].
\end{eqnarray*}

\subsubsection{Validity of the separatrix 
asymptotic expansion}
\par
The separatrix asymptotic expansion 
(\ref{separatrixAsymptoticExpansion}) is valid until:
$$
{\varepsilon U_{n+1}\over U_n }\ll1.
$$
It yields the bounds for the independent 
variable $t$:
$$
\ve\exp(2t)\ll1,\quad 
|t|\ll-{1\over2}\ln(\ve)
$$
and the parameters of the solution $A^+_n$ 
and $B^+_n$:
$$
A^+_n\ll\ve,\quad 
B^+_n\ll\ve\quad \forall n\in{N}.
$$
\begin{lemma}\label{lemmaAboutUpperSeparatrixExpansion}
There exists an asymptotic solution for equation
(\ref{pertubedDuffingsOscillator}) in the form
(\ref{separatrixAsymptoticExpansion}), when $|t|\ll -{1\over2}\ln(\varepsilon)$. 
\end{lemma}

\subsection{Saddle asymptotic 
expansion}
\par
The asymptotic expansion has a form:
\begin{equation}
u(t,\ve)=\pm1+
\sum_{n=1}^\infty\ve^{n/2}
u^\pm_n(\tau)\quad \tau=t+\tau_0
\label{saddleAsymptoticExpansion}
\end{equation}
near the saddle points $u=\pm 1$.
\par
We use the sigh '+' for the expansion near $u=1$ and the sign '-' for the expansion near $u=-1$. 
\par
The correction terms are defined by the following equation:
$$
{u^\pm_n}''-4u^\pm_n=f^\pm_n,
$$
where$f^{\pm}_1\equiv0$,
$f^\pm_2=\cos(\omega\tau-\omega\tau_0+\phi_0)\pm6(u^\pm_1)^2$ and
$f^\pm_n$ is a polynomial of order $3$, which is defined
by the correction terms with the indexes $j,k,l$ such that $j+k+l=n$ when
$n\ge3$. In a general case  $f^\pm_n$ is a finite sum of powers of $e^\tau$, sine,
cosine and independent variable $\tau$. A general formula for $n$-th
correction term  has a form:
\begin{equation}
u^{\pm}_n=\alpha^{\pm}_n \exp(-2\tau)+
\beta^{\pm}_n\exp(2\tau)+ w^{\pm}_n(\tau).
\label{n-thCorrectionTermForSaddleExpansion}
\end{equation}
Here $w^\pm_n(\tau)$ has not the terms $C_1e^{2\tau}$ and
$C_2e^{-2\tau}$ for $\forall C_1,C_2=const$ as $\tau\to\pm\infty$.
\par
First and second corrections are:
\begin{eqnarray*}
u^\pm_1=\alpha^\pm_1\exp(-2\tau)+
\beta^\pm_1\exp(2\tau);
\\
u^\pm_2=\alpha^\pm_2\exp(-2\tau)+
\beta^\pm_2\exp(2\tau)+
{(\alpha^\pm_1)^2\over12}\exp(-4\tau)+
{(\beta^\pm_1)^2\over12}\exp(4\tau) -
\\
- {1\over2}\alpha^\pm_1\beta^+_1-
{1\over 4+\omega^2}
\cos(\omega\tau+\phi_0-\omega\tau_0).
\end{eqnarray*}
\par
The $n$-th correction term is estimated for large $\tau$:
$$
u^\pm_n=O(\exp(\pm2n|\tau|)),\quad n>1,\quad 
\tau\to\pm\infty.
$$
\par
The domain of validity for the saddle asymptotic
expansions is defined by the inequality:
$$
{\ve^{1/2} u^{\pm}_{n+1}\over u^{\pm}_n}\ll1,\quad \tau\to\pm\infty.
$$
It yields:
$$
|\tau|\ll-{1\over 4}\ln(\ve).
$$
\begin{lemma}\label{lemmaAboutSaddleExpansion}
There exists an asymptotic solution in form
(\ref{saddleAsymptoticExpansion}) as $\ve\to0$, where $u_n(\tau)$ has
form (\ref{n-thCorrectionTermForSaddleExpansion}) and  
$|\tau|\ll-{1\over4}\ln(\ve)$.
\end{lemma}

\subsubsection{Matching of the 
asymptotic expansions}
\par
The matching of asymptotic expansions $U(t,\varepsilon)$ and
$u(\tau,\varepsilon)$ yields:
$$
\alpha^+_{2n}=0,\quad \beta^+_{2n}
=0,
$$
\par
One obtains recurrent formulas for the correction terms with odd indexes:
\begin{eqnarray*}
-2\exp(-2t)=\ve^{1/2}\alpha^+_1
\exp(-2\tau), \quad \alpha^+_1=-2,\\
\tau=t+\tau_0,\quad \tau_0=
{1\over4}\ln(\ve).
\end{eqnarray*}
Then:
$$
{1\over16}B_1^+\exp(2t)=
\beta^+_1\ve^{-1/2}\exp(2\tau),
\quad \beta^+_1={1\over16}B_1^+.
$$
For higher-order terms as $t\to\infty$ one obtains:
\begin{eqnarray*}
4A_n^+\exp(-2t)=\ve^{1/2}
\alpha^+_{2n+1}\exp(-2\tau),\\ 
\alpha_{2n+1}^+=4 A_n^+;\\
{1\over16}B_n^+\exp(2t)=
\beta_{2n-1}^+
\ve^{-1/2}\exp(2\tau),\\
\beta_{2n-1}^+={1\over16}B_n^+,
\quad n\in{\mathbf N}.
\end{eqnarray*}
\par
The sign of $\beta_1^+$ depends on the parameter $\phi_0$:
$$
\beta_1^+={1\over16}B_1^+ =
B^-_1 +{1\over16}\cos(\phi_0) 
{1\over\cosh(\pi\omega/2)}
$$
If $\beta_1^+>0$ or the same  $\cos(\phi_0)>-16B^-_1\cosh(\pi\omega/2)$ then the asymptotic
solution of (\ref{pertubedDuffingsOscillator}) goes to infinity,
otherwise when $\cos(\phi_0)<-16B^-_1\cosh(\pi\omega/2)$ the  solution
of (\ref{pertubedDuffingsOscillator}) is close to the separatrix which
goes to the left saddle.
\begin{theorem}
\label{theoremAboutUpperAsymptoticExpansions}
There exists an asymptotic 
solution of (\ref{pertubedDuffingsOscillator}), 
such that this solution has a form 
(\ref{separatrixAsymptoticExpansion}) 
when ${1\over2}\ln(\ve)\ll t\ll
-{1\over2}\ln(\ve)$ and has form 
(\ref{saddleAsymptoticExpansion}) when 
$-{1\over4}\ln(\ve)\ll t\ll -{1\over2}\ln(\ve).$
\end{theorem}

\subsection{Lower separatrix branch}
\par
Let us consider the case when $\beta_1^+<0$, then the  main term of
an asymptotic solution as $\tau\to\infty$ is the separatrix which goes 
from one saddle $(1,0)$ on the plane $(u,u')$ to another saddle
$(-1,0)$. This separatrix is:
$$
u_0=-\tanh(t).
$$
The asymptotic expansion, which is close to the lower separatrix has a
similar form as the asymptotic expansion close to the upper  separatrix:
\begin{equation}
u(\theta,\ve)=u_0(\theta)+
\sum_{n=0}^\infty\ve^n u_n(\theta).
\label{lowerSeparatrixAsymptotics}
\end{equation}
\par
It is easy to see that the asymptotic solutions
(\ref{lowerSeparatrixAsymptotics}) has the same form as
(\ref{separtrixExpansion}), where one should use 
the variable $\theta$ instead of $t$. 
\par
Equation for $n$-th correction term has a following form:
\begin{equation}
u_n''+3u_n-6 u_0^2 u_n=f_n.
\label{linearizedEquation}
\end{equation}
Here $f_n$ is a polynomial of $3$-d order with respect to terms
$u_j,u_k,u_l$, such that  $j+k+l=n$, and 
\begin{equation}
f_n\sim \sum_{k=-\infty}^{n-1} e^{2k\theta}\bigg(\sum_{l=0}^{2(n-1)}
  \theta^l\bigg[\sum_{m=0}^n\bigg(F^{\pm}_{k,l,m}\cos(m\omega
  \theta)+H^\pm_{k,l,m}\sin(m\omega \theta)\bigg)\bigg]\bigg),\quad
  \theta\to\pm\infty,
\label{asymptotticsOfRightHandSide}
\end{equation}
The general formula for the solution of (\ref{linearizedEquation})
yields:
$$
u_n(\theta)\sim \sum_{k=-\infty}^{n} e^{2k\theta}\bigg(\sum_{l=0}^{2n}
  \theta^l\bigg[\sum_{m=0}^{n+1}\bigg(\tilde u^{\pm}_{k,l,m}\cos(m\omega
  \theta)+\tilde v^\pm_{k,l,m}\sin(m\omega \theta)\bigg)\bigg]\bigg),\quad \theta\to\pm\infty.
$$
Denote $a_n^{\pm}=\tilde u_{\mp1,0,0}^{\pm}/4$ and $b_n^{\pm}=16\tilde
u_{\pm1,0,0}^{\pm}$. These parameters define the parameters of
the solution for the large time. The changes of these parameters are:
$$
\Delta a_n=a_n^- -a_n^+,\quad \Delta b_n=b_n^- - b_n^+.
$$

\subsubsection{Validity of separatrix 
expansion}
\par
The lower separatrix expansion is valid 
until:
$$
{\varepsilon u_{n+1}\over u_n }\ll1.
$$
It yields:
$$
\ve\exp(2\theta)\ll1,\quad 
|\theta|\ll-{1\over2}\ln(\ve).
$$

\begin{lemma}
There exists an asymptotic 
solution of 
(\ref{pertubedDuffingsOscillator}) 
in form 
(\ref{lowerSeparatrixAsymptotics}) 
as $|\theta|\ll-{1\over2}\ln(\ve)$.
\end{lemma}

\subsubsection{Matching of lower 
separatrix asymptotic expansion and 
asymptotic expansion close 
to right saddle}
\par
Let $\beta^+_1<0$. In this case one can match the asymptotic 
expansion in the right saddle point and the 
lower separatrix asymptotic expansion 
(\ref{lowerSeparatrixAsymptotics}).
\par
The matching yields:
$$
\ve^{1/2}\beta_1^+\exp(2\tau)=
-2\exp(2\theta) \quad \hbox{and}\quad 
\ve^{1/2}\alpha_1^+\exp(-2\tau)=
\ve 4 a_1^+\exp(-2\theta),
$$
as  
$\tau\to\infty$ and $\theta \to-\infty$. 
Here $\beta_1^+<0$, as a result we 
obtain:
$$
2\tau+{1\over2}\ln(\ve)+
\ln(-\beta_1^+)=
2\theta+\ln2.
$$
Or the same
$$
\theta=\tau+{1\over4}\ln(\ve)+
{1\over2}\ln(-\beta_1^+)-{1\over2}\ln2.
$$
A substitution of  $\theta$ and reductions give a following 
formula:
$$
{1\over2}\ln(\ve)+\ln(\alpha_1^+)-
2\tau= 
\ln(\ve)+2\ln(2)+\ln(a_1^+)-2\tau
-{1\over2}\ln(\ve)-\ln(-\beta_1^+)-
\ln(2)
,
$$
or 
$$
\alpha_1^+(-\beta_1^+)=8a_1^+,
\quad a^+_1={1\over8}\alpha^+_1 
(-\beta^+_1).
$$
It is easy to see that  $\alpha_1^+=-2$. 
A following computation gives formulas 
for $\beta_3^+$  and  $b_1^+$:
$$
\ve^{3/2}\beta_3^+\exp(2\tau)=
{1\over16}\ve b_1^+\exp(2\theta).
$$
As a result we obtain:
$$
\ln(\beta_3^+)+2\tau=-4\ln(2)+
\ln(b_1^+)+2\tau+\ln(-\beta_1^+)-
\ln(2),
$$
or
$$
\beta_3^+={1\over32}b_1^+
(-\beta_1^+),\quad b^+_1=
32{-\beta^+_3\over\beta^+_1}.
$$
The same way leads us to formulas for the 
higher order terms:
\begin{eqnarray*}
\ve^{(2n+1)/2}\alpha_{2n-1}^+
\exp(-2\tau)=4\ve^n a_n^+
\exp(-2\theta),\\
\alpha_{2n-1}^+(-\beta_1^+)=
8 a_n^+,\quad a^+_n=
{1\over8}\alpha^+_{2n-1}
(-\beta^+_1).
\end{eqnarray*}
\begin{eqnarray*}
\ve^{(2n+1)/2}\beta_{2n+1}^+
\exp(2\tau)={1\over16}\ve b_n^+
\exp(2\theta),\\
{\beta_{2n+1}^+\over(-\beta_1^+)}=
{1\over32} b_n^+,\quad b^+_n=
-32{\beta^+_{2n+1}\over\beta^+_1}.
\end{eqnarray*}
\par
It is convenient to write out the formulas which connect the parameters
of the asymptotic expansions near the upper and the lower
separatrices:
\begin{eqnarray*}
a_1^+={1\over16}(B_1^-+\Delta B_1),\\
a_{n+1}^+=-{1\over32}(A_{n}^++\Delta A_{n})(B_1^-+\Delta B_1),\\
b_n^+=-32{B_{n+1}^-+\Delta B_{n+1}\over B_1^-+\Delta B_1},\quad
  n\in{\mathbf N};\\
\theta=t+{1\over2}\ln(\varepsilon)+{1\over2}\ln\bigg(-{1\over16}(B_1^-+\Delta
B_1)\bigg)-{1\over2}\ln(2),\\
\phi=-{\omega\over2}\bigg[\ln(\varepsilon)+\ln\bigg(-{1\over16}(B_1^-+\Delta B_1)\bigg)-\ln(2)\bigg]+\Phi.
\end{eqnarray*}

\subsection{Neighborhood of left 
saddle point}
\subsubsection{Matching of saddle asymptotic expasion and 
lower separatrix asymptotic expansion}
\par
The matching of the expansions near left saddle point yields:
$$
2\exp(-2\theta)=\ve^{1/2}\alpha_1^-
\exp(-2\sigma),
$$
$$
\ve {1\over16}b_1^-\exp(2\theta)=
\ve^{1/2}\beta_1^-\exp(2\sigma),
$$
$$
\ve 4 a_1^-\exp(-2\theta)=\ve^{3/2}
\alpha_3^-\exp(-2\sigma).
$$
As a result one obtain formulas for 
pthe arameters of the saddle asymptotic 
expansion:
$$
-2\theta={1\over2}\ln(\ve)-2\sigma,
$$
$$
\alpha_1^-=2,
$$
$$
4 a_1^-=\alpha_3^-,
$$
$$
\beta_1^-={1\over16}b_1^-.
$$
A matching yields formulas for 
higher-order terms:
$$
\ve^{(2n+1)/2}(\alpha_{2n+1}^-
\exp(-2\sigma)+\beta^-_{2n+1}
\exp(2\sigma))=\ve^n(4a_n^-
\exp(-2\theta)+{1\over16}b_n^-
\exp(2\theta)).
$$
There fore one obtains formulas for 
the coefficients of the asymptotic 
expansion:
$$
\alpha_{2n-1}^-=4a_n^-,\quad 
16\beta_{2n+1}^-=b_n^-.
$$
\par
Here one can obtain an explicit formula 
for the changing of coefficient $b_n$. 
For example,
$$
\Delta b_1=b^-_1-b^+_1={\pi \over 
\cosh(\pi\omega/2)}
\cos\bigg( -{\omega\over2}(\ln(\ve)+
\ln(-\beta^+_1)-\ln(2))+\Phi\bigg).
$$
\par
As a result we obtain the following theorem.
\begin{theorem}
If $\beta_1^+<0$ then there exists 
an asymptotic solution of 
(\ref{pertubedDuffingsOscillator}) when 
${1\over2}\ln(\ve)\ll t\ll-{1\over2}\ln(\ve)$ 
and has following form 
(\ref{separatrixAsymptoticExpansion}) when 
${1\over2}\ln(\ve)\ll t\ll-{1\over2}\ln(\ve)$
(\ref{saddleAsymptoticExpansion}) 
with sub sign "+" when 
$-{1\over4}\ln(\ve)\ll t\ll-{1\over2}\ln(\ve)$,
(\ref{lowerSeparatrixAsymptotics})
$1\ll t\ll-\ln(\ve)$ and 
(\ref{saddleAsymptoticExpansion}) 
with sign "-" when 
$-{3\over4}\ln(\ve)\ll t\ll-\ln(\ve)$.
\end{theorem}
\subsection{Next circle near separatrix}
\par
Let us consider an oscillation near 
the separatrix. Denote parameters of 
asymptotic solution by sequence numbers:
\begin{eqnarray*}
t_0=t,\quad \theta_0=\theta,\quad 
\tau_0=\tau,\\
\sigma_0=\sigma,\quad B^\pm_n(0)=B^\pm_n,\quad
b^\pm_n(0)=b^\pm_n,\\
\Phi(0)=\Phi,\quad \phi(0)=\phi.
\end{eqnarray*}
\par
The parameters of the asymptotic expansion for 
following circles may being calculated 
by formulas:
\par
if $B^-_1(m)+\Delta B_1(m)<0$ and $n\ge 2$:
\begin{eqnarray}
\theta(m)=t(m)+{1\over2}\ln(\varepsilon)+
{1\over2}\ln\bigg({1\over16}(B^-_1(m)+
\Delta B_1(m))\bigg) -{1\over2}\ln(2),\nonumber\\
b^+_n(m)=-32{B^-_{n+1}(m)+\Delta B_{n+1}(m)\over
B^-_1(m)+\Delta B_1(m)}; \label{parametersOfUpperSeparatrix}\\
a^+_1={1\over32}\bigg(B_1^-(m)+\Delta B_1(m)\bigg),\nonumber\\
a^+_n(m)=-{1\over32}\bigg[A_{n-1}^-(m)+\Delta A_{n-1}(m)\bigg]
\bigg(B_1^-(m)+\Delta B_1(m)\bigg);\nonumber
\end{eqnarray}
\par
if $b^+_1(m)+\Delta b_1(m)>0$:
\begin{eqnarray}
t(m+1)=\theta(m)+{1\over2}\ln(\varepsilon)+
\ln\bigg({1\over16}(b^+_1(m)+\Delta b_1(m))\bigg)-
{1\over2}\ln(2),\nonumber\\
B^-_n(m+1)=-32{b^+_{n+1}(m)+\Delta b_{n+1}(m)\over
b^+_1(m)+\Delta b_1(m)},\label{parametersOfLowerSeparatrix}\\
A^-_1={1\over32}\bigg(b_1^+(m)+\Delta b_1(m)\bigg),\nonumber\\
A^-_n(m)=-{1\over32}\bigg[a_{n-1}^+(m)+\Delta a_{n-1}(m)\bigg]
\bigg(b_1^+(m)+\Delta b_1(m)\bigg);\nonumber
\end{eqnarray}
\par
In particular,
\begin{eqnarray}
B^+_1(m)=B^-_1(m)+
{\pi\over16}{\cos(\Phi(m))\over 
\cosh({\pi\omega\over 2})}\nonumber\\
b^-_1(m)=b^+_1(m)+
{\pi\over16}{\cos(\phi(m))\over 
\cosh({\pi\omega\over 2})}.\label{m-thMainTerms}
\end{eqnarray}
\par
For given parameter $\Phi(0)$ we obtain following 
discreet dynamics over a near-separatrix circle:
\begin{eqnarray}
\Phi(m+1)&=\omega(t(m+1)-
\theta(m+1))+\phi(m)\nonumber\\
&=-
{\omega\over2}\bigg(\ln(\varepsilon)+\ln({1\over16}(b^+_1(m)+
\Delta b_1(m))-\ln(2)\bigg) +\phi(m),
\label{m-thSift}
\\
\phi(m)&=\omega(\theta(m)-t(m))+\Phi(m)=
-\nonumber\\
&\omega(\ln(\varepsilon)+
{1\over2}\ln({1\over16}(B^-_1(m)+
\Delta B_1(m))-\nonumber\\
&{1\over2}\ln(2)) +
\Phi(m).\nonumber
\end{eqnarray}
\par
It is easy to see that $\Phi(m)$ and 
$\phi(m)$ have the same order:
\begin{eqnarray*}
\phi(m)=-\omega\bigg({1\over2}+m\bigg)
\ln(\varepsilon)+O(1),\\
\Phi(m)=-m\omega\ln(\varepsilon)+O(1).
\end{eqnarray*}
\par
The shown formulas give a 
dependency on the sequence 
$\{B^-_n\}_{n=0}^\infty$ for 
dynamics of the solution for times 
$t=O(N\ln(\varepsilon))$ for 
$\forall N\in\mathbf N$. 
\par
\begin{lemma}
For $\forall \varepsilon\in(0,\varepsilon_0)$,
where $\varepsilon=\const>0$ there exists 
a sequence $\{B^-_n\}_{n=0}^\infty$ such that 
$B_1^+(k)<0$ and $b_1^+>0$ $\forall k$. 
\end{lemma}
\par
The sequence of the lemma and inequalities 
(\ref{parametersOfUpperSeparatrix}) and 
(\ref{parametersOfLowerSeparatrix}) give
\begin{theorem}
For $\forall \varepsilon\in(0,\varepsilon_0)$,
where $\varepsilon=\const>0$ 
there exists the sequence  
$\{B^-_n\}_{n=0}^\infty$ such that the asymptotic 
solution oscillates near the separatrices as 
$t=O(N\ln(\varepsilon))$ for all 
$\forall N\in\mathbf N$.
\end{theorem}

\section{Discreet dynamical system}
\label{sectionDescreetDynamicalSystem}
\par
Let us consider the discreet dynamical systems which defines by
separatrix map (\ref{parametersOfUpperSeparatrix}),
(\ref{parametersOfLowerSeparatrix}).  It is convenient to change this
system by the following way. Let us define
\begin{eqnarray*}
\sigma_k(2n-1)\equiv B^-_k(n-1),\quad \sigma_k(2n)\equiv b^+_k(n-1),\quad n\in{\mathbf N};\\
\Delta\sigma_k(2n-1)\equiv\Delta B_k(n),\quad \Delta\sigma_k(2n)\equiv\Delta b_k(n-1);\\
\chi_k(2n-1)\equiv A^-_k(n-1),\chi_k(2n)\equiv a^+_k(n-1),\quad k\ge2,\quad n\in{\mathbf N};\\
\Delta\chi_k(2n-1)\equiv\Delta A_k(n-1),\quad \Delta\chi_k(2n)\equiv\Delta a_k(n-1);\\
\psi(2n-1)\equiv\Phi(n)+n\omega\ln(\varepsilon),\quad \psi(2n)\equiv\phi(n)+(n+{1\over2})\omega\ln(\varepsilon).
\end{eqnarray*}
\par
This formulas is prolonged on the next step if
$(\sigma_1(n)+\Delta\sigma_1(m))(-1)^n>0$. The formulas for a discreet dynamical system are:
\begin{eqnarray*}
\sigma_k(1)=B_k^-,\quad \chi_k(1)=A_k^-,\quad \psi(1)=\Phi,\\
\hbox{if}\,\,(\sigma_1(m)+\Delta\sigma_1(m))(-1)^m>0,\quad \forall
m\le n,\,\,\hbox{then}\\
\sigma_k(n+1)=-32{\sigma_{k+1}(n)+\Delta\sigma_{k+1}(n)\over \sigma_1(n)+\Delta\sigma_1(n)},\\
\psi(n+1)=-{\omega\over2}\bigg(\ln\bigg({1\over16}(\sigma_1(n)+\Delta\sigma_1(n))\bigg)-\ln(2)\bigg)+\psi(n),\\
\chi_1(2n)=2,\quad \chi_1(2n-1)=-2,\\
\chi_k(n+1)=-{1\over32}\bigg(\chi_{k-1}(n)+\Delta\chi_{k-1}(n)\bigg)\big(\sigma_1(n)+\Delta\sigma_1(n)\big).
\end{eqnarray*}
\par
Those formulas give a nonlinear discreet 
dynamical system. The non-linearity defines by 
the terms $\Delta\sigma_k(n),\Delta\chi_k(n)$ or the same 
$\Delta B_k(n)$, $\Delta b_k(n)$, $\Delta A_k(n)$, $\Delta a_k(n)$.
 
\subsection{Generalized Bernoulli shift}
\par 
The separatrix map defines a Bernoulli 
shift for parameters of this map. Formulas 
(\ref{parametersOfUpperSeparatrix}) and 
(\ref{parametersOfLowerSeparatrix}) show that 
the $\sigma_k(n+1)$ depends on  $\sigma_{k+2}(n)$. 
It means that $(n)$-th correction 
on the $(m+1)$-th oscillation depends on 
$(n+2)$-th correction on the $m$-th 
oscillation near the separatrix. This is typical 
property of the Bernoulli shift. 
This property shows the loss of the accuracy for 
the approximation of the motion and the instability 
the solution with respect to initial data.

\subsection{Cantor manifold}
When we study the asymptotic solution we say that there exists
$\varepsilon_0>0$ such that for some $\varepsilon\in(0,\varepsilon_0)$
there exists studying asymptotic solution. In  this subsection we turn
the  study by the other side. We concentrate on the problem for
the structure of the manifold $\varepsilon\in(0,\varepsilon_0)$ which gives
the oscillations near the separatrices. 
\par
The following condition 
\begin{equation}
\cos(\psi(n)-\omega{n-1\over2}\ln(\varepsilon))<-16 \sigma_1(n)\cosh\bigg({\pi\omega\over2}\bigg)\label{condOfProl}
\end{equation}
defines the possibility to prolong of the discreet dynamical system
on $(n+1)$-th step. Define $\delta=-\omega\ln(\varepsilon)$, then
$\delta\in (\delta_0,\infty)$ where
$\delta_0=-\omega\ln(\varepsilon_0)$ and $0<\varepsilon_0<1$.  Let
$\sigma_1(n)$ is such that 
$$
\arccos(-16\sigma_1(n)\cosh({\pi\omega\over2}))-\psi(n)=\pm{\pi\over3}+2\pi
k.
$$
Then (\ref{condOfProl}) defines the parameter $\delta$  
$$
\forall k>\big[ {\omega(n-2)\over4\pi}\ln(\varepsilon_0)\big]+1.
$$
Let 
$$
{n-1\over2}\delta\in(-\pi+2\pi k,\pi+2\pi k)\in(\delta_0,\infty),
$$
then when 
$$
{n-1\over2}\delta<-{\pi\over3}+2\pi k,\,\,\hbox{and}\,\,{n-1\over2}\delta>{\pi\over3}+2\pi k
$$
the discrete dynamical system is prolonged by the next step.
\par
Let us consider sequences $\alpha\in{\mathbf N}$ such that 
$$
-16\sigma(l)\cosh({\pi\omega\over2})>1,\quad \hbox{as}\quad l>n_\alpha,\quad l<3n_\alpha-2
$$
and $n_{\alpha+1}=3 n_{\alpha}-2$,
$$
\arccos(-16\sigma_1(n_{\alpha+1})\cosh({\pi\omega\over2}))-\psi(n_{\alpha+1})=\pm{\pi\over3}+2\pi
k,
$$
$$
\forall
k>\big[{\omega(n_{\alpha+1}-1)\over4\pi}\ln(\varepsilon_0)\big]+1.
$$
Let $\varepsilon\in(0,\varepsilon_0)$, then the set of $\delta$, for
which the discrete dynamical system is prolonged, is the Cantor
set on any $(-\pi+2\pi, k\pi+2\pi k)$.


\begin{thebibliography}{cc}
\bibitem{Poincare}
H. Poincare, Les methods nouvellesde 
la mecanicue celeste,3, Gauthier-Villars, 
Paris, 1899. 
\bibitem{Melnikov}
V.K. Mel'nikov, On the stability of the 
center for time periodic perturbations. 
Trans. Moscow Math. Soc., 1963, v.12, 
pp.1-57.
\bibitem{FilonenkoSagdeevZaslavskii}
N.N. Filonenko, R.Z. Sagdeev, G.M. Zaslavskii, 
Nuclear Fussion, 1967, v.7, p.253.
\bibitem{Neishtadt}
A.I. Neishtadt, Passage trough a 
separatrix in a resonance problem 
with a slowly-varying parameter. 
J. Appl. Math. Mech., 1975, v.39, 
pp. 594-605.
\bibitem{Timofeev}
A.V. Timofeev, On the constancy 
of an adiabatic invariant when the 
nature of motion changes. JEPTH, 
1978, v.48, pp.656-659.
\bibitem{CaryScodjie}
J.R. Cary, R.T. Scodje, Phase change 
between separatrix crossing. Physica D,
v.36, pp.287-316.
\bibitem{DiminieHaberman1}
D.C. Diminnie, R. Haberman, Slow 
passage through a saddle-center 
bifurcation. J. Nonlinear Sci., 2000, v.10, 
pp. 197-221. 
\bibitem{Kiselev}
O.M.Kiselev. Hard Loss of Stability in 
Painleve-2 Equation. Journal of Nonlinear 
Mathematical Physics, 2001, v.8, n1, p.65-95.
\bibitem{KiselevGlebov}
O.M. Kiselev and S.G. Glebov, An 
asymptotic solution slowly crossing 
the separatrix near a saddle-centre 
bifurcation point. Nonlinearity, 2003, 
v.16, pp.327-362. 
\end{thebibliography}
\end{document}